\documentclass{amsart}

\input{epsf.tex}
\usepackage{amsfonts,amssymb,verbatim,amsmath,amsthm,latexsym,textcomp,amscd}
\usepackage{latexsym,amsfonts,amssymb,epsfig,verbatim}
\usepackage{amsmath,amsthm,amssymb,latexsym,graphics,textcomp}
\usepackage{eucal,eufrak}
\usepackage{graphicx}
\usepackage{color}
\usepackage{url}

\input xy
\xyoption{all}

\begin{document}

\newtheorem{theorem}{Theorem}[section]
\newtheorem{prop}[theorem]{Proposition}
\newtheorem{lemma}[theorem]{Lemma}
\newtheorem{cor}[theorem]{Corollary}
\newtheorem{definition}[theorem]{Definition}
\newtheorem{conj}[theorem]{Conjecture}
\newtheorem{rmk}[theorem]{Remark}
\newtheorem{qn}[theorem]{Question}
\newtheorem{claim}[theorem]{Claim}
\newtheorem{defth}[theorem]{Definition-Theorem}

\newcommand{\boundary}{\partial}
\newcommand{\C}{{\mathbb C}}
\newcommand{\integers}{{\mathbb Z}}
\newcommand{\natls}{{\mathbb N}}
\newcommand{\ratls}{{\mathbb Q}}
\newcommand{\reals}{{\mathbb R}}
\newcommand{\proj}{{\mathbb P}}
\newcommand{\lhp}{{\mathbb L}}
\newcommand{\tube}{{\mathbb T}}
\newcommand{\cusp}{{\mathbb P}}
\newcommand\AAA{{\mathcal A}}
\newcommand\BB{{\mathcal B}}
\newcommand\CC{{\mathcal C}}
\newcommand\DD{{\mathcal D}}
\newcommand\EE{{\mathcal E}}
\newcommand\FF{{\mathcal F}}
\newcommand\GG{{\mathcal G}}
\newcommand\HH{{\mathcal H}}
\newcommand\II{{\mathcal I}}
\newcommand\JJ{{\mathcal J}}
\newcommand\KK{{\mathcal K}}
\newcommand\LL{{\mathcal L}}
\newcommand\MM{{\mathcal M}}
\newcommand\NN{{\mathcal N}}
\newcommand\OO{{\mathcal O}}
\newcommand\PP{{\mathcal P}}
\newcommand\QQ{{\mathcal Q}}
\newcommand\RR{{\mathcal R}}
\newcommand\SSS{{\mathcal S}}
\newcommand\TT{{\mathcal T}}
\newcommand\UU{{\mathcal U}}
\newcommand\VV{{\mathcal V}}
\newcommand\WW{{\mathcal W}}
\newcommand\XX{{\mathcal X}}
\newcommand\YY{{\mathcal Y}}
\newcommand\ZZ{{\mathcal Z}}
\newcommand\CH{{\CC\HH}}
\newcommand\PEY{{\PP\EE\YY}}
\newcommand\MF{{\MM\FF}}
\newcommand\RCT{{{\mathcal R}_{CT}}}
\newcommand\RCTT{{{\mathcal R}^2_{CT}}}
\newcommand\PMF{{\PP\kern-2pt\MM\FF}}
\newcommand\FL{{\FF\LL}}
\newcommand\PML{{\PP\kern-2pt\MM\LL}}
\newcommand\GL{{\GG\LL}}
\newcommand\Pol{{\mathcal P}}
\newcommand\half{{\textstyle{\frac12}}}
\newcommand\Half{{\frac12}}
\newcommand\Mod{\operatorname{Mod}}
\newcommand\Area{\operatorname{Area}}
\newcommand\ep{\epsilon}
\newcommand\hhat{\widehat}
\newcommand\Proj{{\mathbf P}}
\newcommand\U{{\mathbf U}}
 \newcommand\Hyp{{\mathbf H}}
\newcommand\D{{\mathbf D}}
\newcommand\Z{{\mathbb Z}}
\newcommand\R{{\mathbb R}}
\newcommand\Q{{\mathbb Q}}
\newcommand\E{{\mathbb E}}
\newcommand\til{\widetilde}
\newcommand\length{\operatorname{length}}
\newcommand\tr{\operatorname{tr}}
\newcommand\gesim{\succ}
\newcommand\lesim{\prec}
\newcommand\simle{\lesim}
\newcommand\simge{\gesim}
\newcommand{\simmult}{\asymp}
\newcommand{\simadd}{\mathrel{\overset{\text{\tiny $+$}}{\sim}}}
\newcommand{\ssm}{\setminus}
\newcommand{\diam}{\operatorname{diam}}
\newcommand{\pair}[1]{\langle #1\rangle}
\newcommand{\T}{{\mathbf T}}
\newcommand{\inj}{\operatorname{inj}}
\newcommand{\pleat}{\operatorname{\mathbf{pleat}}}
\newcommand{\short}{\operatorname{\mathbf{short}}}
\newcommand{\vertices}{\operatorname{vert}}
\newcommand{\collar}{\operatorname{\mathbf{collar}}}
\newcommand{\bcollar}{\operatorname{\overline{\mathbf{collar}}}}
\newcommand{\I}{{\mathbf I}}
\newcommand{\tprec}{\prec_t}
\newcommand{\fprec}{\prec_f}
\newcommand{\bprec}{\prec_b}
\newcommand{\pprec}{\prec_p}
\newcommand{\ppreceq}{\preceq_p}
\newcommand{\sprec}{\prec_s}
\newcommand{\cpreceq}{\preceq_c}
\newcommand{\cprec}{\prec_c}
\newcommand{\topprec}{\prec_{\rm top}}
\newcommand{\Topprec}{\prec_{\rm TOP}}
\newcommand{\fsub}{\mathrel{\scriptstyle\searrow}}
\newcommand{\bsub}{\mathrel{\scriptstyle\swarrow}}
\newcommand{\fsubd}{\mathrel{{\scriptstyle\searrow}\kern-1ex^d\kern0.5ex}}
\newcommand{\bsubd}{\mathrel{{\scriptstyle\swarrow}\kern-1.6ex^d\kern0.8ex}}
\newcommand{\fsubeq}{\mathrel{\raise-.7ex\hbox{$\overset{\searrow}{=}$}}}
\newcommand{\bsubeq}{\mathrel{\raise-.7ex\hbox{$\overset{\swarrow}{=}$}}}
\newcommand{\tw}{\operatorname{tw}}
\newcommand{\base}{\operatorname{base}}
\newcommand{\trans}{\operatorname{trans}}
\newcommand{\rest}{|_}
\newcommand{\bbar}{\overline}
\newcommand{\UML}{\operatorname{\UU\MM\LL}}
\newcommand{\EL}{\mathcal{EL}}
\newcommand{\tsum}{\sideset{}{'}\sum}
\newcommand{\tsh}[1]{\left\{\kern-.9ex\left\{#1\right\}\kern-.9ex\right\}}
\newcommand{\Tsh}[2]{\tsh{#2}_{#1}}
\newcommand{\qeq}{\mathrel{\approx}}
\newcommand{\Qeq}[1]{\mathrel{\approx_{#1}}}
\newcommand{\qle}{\lesssim}
\newcommand{\Qle}[1]{\mathrel{\lesssim_{#1}}}
\newcommand{\simp}{\operatorname{simp}}
\newcommand{\vsucc}{\operatorname{succ}}
\newcommand{\vpred}{\operatorname{pred}}
\newcommand\fhalf[1]{\overrightarrow {#1}}
\newcommand\bhalf[1]{\overleftarrow {#1}}
\newcommand\sleft{_{\text{left}}}
\newcommand\sright{_{\text{right}}}
\newcommand\sbtop{_{\text{top}}}
\newcommand\sbot{_{\text{bot}}}
\newcommand\sll{_{\mathbf l}}
\newcommand\srr{_{\mathbf r}}
\newcommand\geod{\operatorname{\mathbf g}}
\newcommand\mtorus[1]{\boundary U(#1)}
\newcommand\A{\mathbf A}
\newcommand\Aleft[1]{\A\sleft(#1)}
\newcommand\Aright[1]{\A\sright(#1)}
\newcommand\Atop[1]{\A\sbtop(#1)}
\newcommand\Abot[1]{\A\sbot(#1)}
\newcommand\boundvert{{\boundary_{||}}}
\newcommand\storus[1]{U(#1)}
\newcommand\Momega{\omega_M}
\newcommand\nomega{\omega_\nu}
\newcommand\twist{\operatorname{tw}}
\newcommand\modl{M_\nu}
\newcommand\MT{{\mathbb T}}
\newcommand\Teich{{\mathcal T}}
\renewcommand{\Re}{\operatorname{Re}}
\renewcommand{\Im}{\operatorname{Im}}

\title{On Discreteness of Commensurators}

\author{Mahan Mj}

\address{RKM Vivekananda University, Belur Math, WB-711 202, India}

\date{}

\begin{abstract} 
We begin by showing that commensurators of  Zariski dense subgroups of isometry groups of  symmetric spaces of non-compact type 
are  discrete provided that the limit set on the Furstenberg boundary is not invariant under the action of a (virtual) simple factor.
In particular for rank one or simple Lie groups, Zariski dense subgroups with non-empty domain of discontinuity have discrete commensurators.
This generalizes a Theorem of Greenberg for Kleinian groups.
 We then prove that for  all finitely generated, Zariski dense, infinite covolume discrete subgroups of 
$Isom ({\mathbb{H}}^3)$, commensurators are  discrete. Together these prove discreteness of commensurators for all known examples of finitely
generated, Zariski dense, infinite covolume discrete subgroups of $Isom(X)$ for $X$ a symmetric space of non-compact type.

\end{abstract}

\maketitle

\begin{center}
AMS subject classification =   57M50

\end{center}

\tableofcontents

\section{Introduction}
Motivated by  Margulis'  celebrated   characterization of arithmeticity of irreducible lattices in semi-simple Lie groups
in terms of density of the commensurator \cite{margulis}, Shalom asked (cf. \cite{llr})for a description of commensurators of
 Zariski dense subgroups of  semi-simple Lie groups.
We start with the observation that  for Zariski dense subgroups of isometry groups of
 rank one symmetric spaces, commensurators are discrete provided that the domain of discontinuity is non-empty.

\medskip

\noindent {\bf Proposition \ref{zariski} and Corollary \ref{zariskicor}:}   {\it Suppose $\Gamma$ is a  Zariski-dense subgroup
of a    semi-simple Lie group $L = Isom(X)$ for $X$ a rank one symmetric space
of non-compact type.  Then the commensurator
$Comm(\Gamma)$ of $\Gamma$  is  discrete if the limit set $\Lambda_\Gamma$  is not all 
of $\partial X$.}

\medskip

When $X = {\mathbb H}^3$,  the above is due to Greenberg \cite{greenberg1}, \cite{greenberg2}.

Using the theory of limit sets (in the Furstenberg boundary $Y = G/B$)
of Zariski-dense subgroups of semi-simple Lie groups of higher rank developed by   Benoist \cite{benoist-gafa}
we generalize Theorem \ref{zariski} to arbitrary semi-simple Lie groups.

\medskip

{\bf Proposition \ref{zariskicor2}:} {\it Suppose $\Gamma$ is a Zariski-dense subgroup
of a  semi-simple Lie group $L = Isom(X)$ for $X$ a  symmetric space
of non-compact type.  Further suppose that the limit set $\Lambda_\Gamma \subset Y$  is not invariant under any non-trivial semi-simple
(virtual) factor
 subgroup $L_0$. Then the commensurator
$Comm(\Gamma)$ of $\Gamma$  is   discrete  in $L$.   }

\medskip

We then specialize to the case that $\Gamma$ is abstractly a (relative) hyperbolic group in the sense of Gromov and $X$ is of rank one.

Infinite covolume discrete
subgroups of semi-simple Lie groups have received  little attention in the context
of commensurators until  recently. However, the theory of
Kleinian groups deals primarily with infinite covolume discrete
subgroups of $PSL_2({\mathbb C})$.
Commensurators of some non-free Kleinian groups have been investigated by Leininger, Long and Reid
\cite{llr}. In this paper we reprove and extend
 their results to all finitely generated Kleinian groups.

\medskip

\noindent {\bf Theorem \ref{comm2prime}:}    {\it
Let $G$ be  a finitely generated, Zariski dense Kleinian group such that
  ${\mathbb{H}}^3/G$ has infinite volume;
then the commensurator $Comm(G)$ of $G$ is discrete in $PSL_2({\mathbb{C}})$. }

\medskip

The above statement 
 was proven by Leininger, Long and Reid \cite{llr} when $G$ is
non-free and without parabolics. As in \cite{llr} Theorem \ref{comm2prime} can be strengthened as
follows.

\medskip

\noindent {\bf Theorem \ref{comm3}:}{\it Let $G$ be a finitely generated Kleinian group such that
 $G$ is Zariski dense in $PSL_2({\mathbb{C}})$ and ${\mathbb{H}}^3/G$ has infinite volume.
Then  $[Comm(G):G] < \infty$ unless $G$ is virtually a fiber subgroup, in which case
$Comm(G)$ is the fundamental group of a virtually fibered finite volume hyperbolic 3-manifold. }

\medskip

Together, Proposition \ref{zariskicor2} and  Theorem \ref{comm2prime} prove discreteness of commensurators for all known examples of finitely
generated, Zariski dense, infinite covolume discrete subgroups of $Isom(X)$ for $X$ a symmetric space of non-compact type. 

We also provide an example of a Zariski dense subgroup $G \subset Isom ({\mathbb{H}}^4)$ such that
$Comm(G)$ is discrete, $[Comm(G):G] = \infty$
but the limit set $\Lambda_G \neq \partial {\mathbb{H}}^4$. We also justify a statement of Gromov \cite{gromov-sq} by showing that
negatively curved 4-manifolds cannot fiber over the circle.


\subsection{Zariski-dense subgroups of semi-simple Lie groups} Our first observation is that a version of Margulis' dichotomy
holds for all Zariski dense subgroups. (Some version of this is probably well-known to experts in the theory of Lie groups and Algebraic Groups.)

\begin{prop} Suppose $\Gamma$ is a Zariski-dense subgroup
of a  semi-simple Lie group $L = Isom(X)$ for $X$ a  rank one symmetric space
of non-compact type.  Then the commensurator
$Comm(\Gamma)$ of $\Gamma$  is either  discrete or dense in $L$. \label{zariski} \end{prop}

\noindent{\bf Proof:} Let $ \overline{Comm(\Gamma)}$ be the closure of $Comm(\Gamma)$
in the Lie group $L$. If the connected component of the identity $L_0$ of $ \overline{Comm(\Gamma)}$ is trivial,
then $ \overline{Comm(\Gamma)}$ and hence ${Comm(\Gamma)}$ is totally disconnected and therefore discrete.

Else, let ${\mathfrak{l}}_0$ denote the Lie algebra of $L_0$ and 
$\mathfrak{l}$ denote the Lie algebra of $L$. Since $L_0$ is the connected component of the identity  of $ \overline{Comm(\Gamma)}$,
$L_0$ is normal in $ \overline{Comm(\Gamma)}$. Then  ${\mathfrak{l}}_0$ is an invariant subspace of $\mathfrak{l}$ under the adjoint
action $Ad_\Gamma$ of $\Gamma$. Hence  ${\mathfrak{l}}_0$ is an invariant subspace of $\mathfrak{l}$ under the adjoint
action of $ZC(\Gamma )$, the Zariski closure of $\Gamma$. (This follows from the fact that $Ad_G$ acts on the Grassmannian of $dim({{\mathfrak{l}}}_0)$ - planes 
in $\mathfrak{l}$ and $Ad_\Gamma$ fixes the plane ${{\mathfrak{l}}}_0$.)  Since $\Gamma$ is Zariski-dense, it follows that 
${{\mathfrak{l}}}_0$ is an ideal in ${{\mathfrak{l}}}$. Let ${{\mathfrak{l}}} = {{\mathfrak{g}}}_1 \oplus {\mathfrak{g}}_2 \oplus \cdots 
\oplus {{\mathfrak{g}}}_k$ be the decomposition of ${{\mathfrak{l}}}$ into its simple Lie algebra summands. Hence
${{\mathfrak{l_0}}} = {{\mathfrak{g}}}_{i_1} \oplus {\mathfrak{g}}_{i_2} \oplus \cdots 
\oplus {{\mathfrak{g}}}_{i_s}$ for $ i_1, \cdots , i_s $ distinct elements of $\{ 1, \cdots , k\}$.

2 cases arise: \\
{\bf Case 1:} $L_0$ is non-compact. Let $K \subset L$ be the maximal compact subgroup of $L$
and $X = L/K$ be the associated symmetric space. Then $X_0 = L_0 / (L_0 \cap K)$ is the symmetric 
subspace of $L$ associated to the subgroup $L_0$. Further, $L_0$ is normal in $L$, since 
${\mathfrak{l}}_0$ is an ideal in ${\mathfrak{l}}$ and $L_0$ is connected. Since $L_0$ is non-compact, its limit set is non-empty and since $L_0$ is normal,
its limit set is all of $\partial X$ (This is essentially the only place where we really use the rank one assumption).

Using the associated Killing form \cite{helgason}
for instance, $X_0$ is a totally geodesic subspace of
$X$. But from the previous paragraph, the limit set of $L_0$ is all of $\partial X$. Hence $X_0 = X$ and $ \overline{Comm(\Gamma)} = L = Isom(X)$. \\

{\bf Case 2:} $L_0$ is compact. Then $L_0$ fixes some point $x \in X$. Further, $L_0$ is normal in $L$, since 
${\mathfrak{l}}_0$ is an ideal in ${\mathfrak{l}}$. Hence for any $g \in L$ and $l \in L_0$, $(glg^{-1}) g.x = gl.x = g.x$ and so
$L_0$ fixes $g.x$. Hence $L_0$ fixes all $x \in X$. Therefore $L_0$ is trivial, since $L= Isom(X)$. $\Box$

\medskip

As an immediate Corollary of Theorem \ref{zariski}, we have the following:

\begin{cor} Suppose $\Gamma$ is a Zariski-dense subgroup
of a  semi-simple Lie group $L = Isom(X)$ for $X$ a rank one symmetric space
of non-compact type.  Further suppose that the limit set $\Lambda_\Gamma \neq \partial X$. Then the commensurator
$Comm(\Gamma)$ of $\Gamma$  is   discrete  in $L$. \label{zariskicor} \end{cor}

\noindent{\bf Proof:} The limit set $\Lambda_\Gamma$ is invariant under $Comm (\Gamma )$. Hence if $Comm (\Gamma )$ is dense in $L$, 
$\Lambda_\Gamma = \partial X$. Therefore by Theorem \ref{zariski}, $Comm(\Gamma)$   is   discrete  in $L$. $\Box$

\medskip

Next, we generalize Proposition \ref{zariski} to arbitrary semi-simple $L= Isom(X)$, by introducing an appropriate generalization of the
notion of irreducibility. Motivated by the definition of irreducibility of lattices we propose the following.

\begin{definition} A Zariski dense subgroup $\Gamma$ of a semi-simple Lie group $L$ will be called {\bf strongly irreducible} if 
for any proper  non-trivial normal subgroup $L_0$ of $L$, the projection of $\Gamma$ to the quotient $L/L_0$ is indiscrete. \end{definition}

Now, let $\Gamma$ be a strongly irreducible Zariski-dense subgroup of a semi-simple Lie group $L = Isom (X)$. As in Theorem \ref{zariski}
let $L_0$ be the connected component of the identity in $\overline{Comm(\Gamma )}$.
Suppose $L_0 \neq \{ 1 \}$.  The quotient $\overline{Comm(\Gamma )} /L_0$
is discrete. Since $\Gamma \subset \overline{Comm(\Gamma )}$, it follows that $\overline{Comm(\Gamma )} = L_0$. Since $\Gamma$ is Zariski-dense,
$L_0 = L$ and we have shown the following.

\begin{prop} Suppose $\Gamma$ is a strongly  irreducible Zariski-dense subgroup
of a  semi-simple Lie group $L = Isom(X)$ for $X$ a   symmetric space
of non-compact type.  Then the commensurator
$Comm(\Gamma)$ of $\Gamma$  is either  discrete or dense in $L$. \label{zariski3} \end{prop}

The following is another simple generalization of Theorem \ref{zariski}.

\begin{prop} Suppose $\Gamma$ is a Zariski-dense subgroup
of a  semi-simple Lie group $L = Isom(X)$ for $X$ a   symmetric space
of non-compact type.  If $Comm(\Gamma)$ denotes the commensurator
 of $\Gamma$, then there exists a short exact sequence
$$1 \rightarrow L_0 \rightarrow \overline{Comm(\Gamma)} \rightarrow Q \rightarrow 1$$
where $L_0$ is a semi-simple Lie subgroup of $L$ and $Q$   is   discrete  in $L/L_0$. \label{zariski4} \end{prop}

\noindent{\bf Proof:} As in the proof of Theorem \ref{zariski},
let $ \overline{Comm(\Gamma)}$ be the closure of $Comm(\Gamma)$
in the Lie group $L$ and $L_0$ be the connected component of the identity in  $ \overline{Comm(\Gamma)}$. If $L_0$ is trivial,
then $ \overline{Comm(\Gamma)}$  is totally disconnected and hence ${Comm(\Gamma)}$ is discrete.

Else, as before, let ${\mathfrak{l}}_0$ denote the Lie algebra of $L_0$ and 
$\mathfrak{l}$ denote the Lie algebra of $L$. Then $L_0 \subset \overline{Comm(\Gamma)}$ is normal
and $ \overline{Comm(\Gamma)} /L_0 = Q$ is discrete. $\Box$

\medskip

We now proceed to refine some of these results to higher rank semi-simple Lie groups in the context
of Zariski dense subgroups of $Isom(X)$ acting on the flag variety or the Furstenberg boundary.
This gives a geometric perspective on Proposition \ref{zariski4}. Let $G$ be a semi-simple Lie group, $K$
a maximal compact subgroup and and $P$ a Borel subgroup. Then $G/P$ is called the
Furstenberg boundary of the associated symmetric space $X = G/K$. See \cite{benoist-gafa}
for details. Limit sets of Zariski-dense subgroups $\Gamma$ of $G$ have been defined by Benoist \cite{benoist-gafa} as in the rank one case.
The limit set of the action of $\Gamma$ on $G/P$  will be denoted as $\Lambda ( G/P, \Gamma)$. The following basic Proposition is due to Benoist,
generalizing an analogous statement for $SL_n({\mathbb{R}})$ by Guivarc'h.

\begin{prop} (Benoist \cite{benoist-gafa}) Let $\Gamma \subset G$ be
 a Zariski dense subgroup. Then  $\Lambda ( G/P, \Gamma)$
is the unique minimal closed $\Gamma$-invariant subset $G/P$. \label{minimal} \end{prop}

A semi-simple Lie algebra is a direct sum of simple Lie algebras. A closed semi-simple Lie subgroup $L_0$ of a Lie group $L$ will be called
a {\bf (virtual) factor} if its Lie algebra ${\mathfrak{l}}_0$ is a direct summand of the Lie algebra  ${\mathfrak{l}}$ of $L$.

Proposition \ref{minimal} will be essential in generalizing Corollary \ref{zariski} to the following.

\begin{prop} Suppose $\Gamma$ is a Zariski-dense subgroup
of a  semi-simple Lie group $L = Isom(X)$ for $X$ a  symmetric space
of non-compact type.  Further suppose that the limit set $\Lambda ( G/P, \Gamma)$ is not invariant under any non-trivial semisimple
(virtual) factor $L_0$. Then the commensurator
$Comm(\Gamma)$ of $\Gamma$  is   discrete  in $L$. \label{zariskicor2} \end{prop}

\noindent{\bf Proof:} As in the proof of Theorem \ref{zariski},
let $ \overline{Comm(\Gamma)}$ be the closure of $Comm(\Gamma)$
in the Lie group $L$ and $L_0$ be the connected component of the identity in  $ \overline{Comm(\Gamma)}$. If $L_0$ is trivial,
then $ \overline{Comm(\Gamma)}$ is totally disconnected and hence ${Comm(\Gamma)}$  is discrete.

Else, as before, let ${\mathfrak{l}}_0$ denote the Lie algebra of $L_0$ and 
$\mathfrak{l}$ denote the Lie algebra of $L$.

2 cases arise: \\
{\bf Case 1:} $L_0$ is compact. Then exactly the same proof as Case 2 of Proposition \ref{zariski} shows that $L_0$ is trivial.\\
{\bf Case 2:} $L_0$ is non-compact. $Comm(\Gamma ) \cap L_0 = \Gamma_0$ is an infinite normal subgroup of $Comm(\Gamma )$.
Also, $\Gamma_0$ is Zariski-dense in $L_0$. Let $\Lambda_\Gamma$ and $\Lambda_{\Gamma_0}$ denote the limit sets
of $\Gamma$ and $\Gamma_0$ respectively. Both are non-empty. 

For every $g \in L$, $g. \Lambda_\Gamma$ is the limit set of $g.\Gamma$ and hence $g.\Gamma g^{-1}$. Since the limit set of a Zariski dense
group is the unique minimal invariant set (Proposition \ref{minimal}), and since a finite index subgroup of a Zariski dense group is Zariski dense,
 $g. \Lambda_\Gamma =  \Lambda_\Gamma$ for every $g \in Comm (\Gamma )$, and hence for every $g \in \overline{Comm(\Gamma)}$. In particular,
$\Lambda_\Gamma$ is $L_0$-invariant. But by hypothesis $\Lambda_\Gamma$ is not invariant under any such $L_0$, a contradiction.

 Therefore  $Comm(\Gamma)$   is   discrete  in $L$. $\Box$

\begin{rmk} In terms of limit sets, invariance of $\Lambda_\Gamma \subset G/P$ under the action of a semi-simple Lie subgroup $L_0 \subset L$
means that $\Lambda_\Gamma = L_0/P_0 \times \Lambda_1$ , where \\
1) $L = L_0 \times L_1$ (at least virtually) \\
2) $L_0/P_0$ is the Furstenberg boundary of $L_0$ \\
3) $\Lambda_1$ is the limit set of the
induced action of $\Gamma$ on the Furstenberg boundary $L_1/P_1$ of the semisimple $L_1$. \end{rmk}



\section{Relations on Boundaries}

\subsection{Hyperbolicity and Relative Hyperbolicity} \label{rh}

\begin{definition}
For any geodesic metric space
$(H,d)$, the {\em hyperbolic cone}
$H^h$ is the metric space
$H\times [0,\infty) = H^h$ equipped with a
path metric $d_h$ obtained from two pieces of
 data \\
1) $d_{h,t}((x,t),(y,t)) = 2^{-t}d_H(x,y)$, where $d_{h,t}$ is the induced path
metric on $H\times \{t\}$.  Paths joining
$(x,t),(y,t)$ and lying on  $H\times \{t\}$
are called {\em horizontal paths}. \\
2) $d_h((x,t),(x,s))=\vert t-s \vert$ for all $x\in H$ and for all $t,s\in [0,\infty)$, and the corresponding paths are called
{\em vertical paths}. \\
3)  for all $x,y \in H^h$,  $d_h(x,y)$ is the path metric induced by the collection of horizontal and vertical paths. \\
\end{definition}

\begin{definition} \label{gro-relhyp}
Let $X$ be a proper (i.e. complete and locally compact) geodesic metric space and $\HH$ be a collection of  uniformly separated subsets of $X$.
$X$ is said to be  hyperbolic relative to $\HH$ in the sense of Gromov, if the  space $\GG (X, \HH)$,  obtained by attaching the hyperbolic cones
$ H^h$ to $H \in \HH$  by identifying $(z,0)$ with $z$
for all $H\in \HH$ and $z \in H$,
 is a proper hyperbolic metric space. The collection $\{ H^h : H \in \HH \}$ is denoted
as ${\HH}^h$. The induced path metric on $\GG (X, \HH)$ is also denoted by $d_h$.
\end{definition}

A group $G$ is  hyperbolic relative to a finite collection $H_1, \cdots , H_k$, if
the Cayley graph
$\Gamma$ of $G$ with respect to some finite generating set
 is  hyperbolic relative to the collection $\HH$ of translates of Cayley subgraphs of $H_1, \cdots , H_k$
(see
\cite{gromov-hypgps},
\cite{farb-relhyp} and \cite{bowditch-relhyp} for details on relative hyperbolicity). Then $\partial G
= \partial \GG(G,\HH )$ is called
the (relative) hyperbolic (or Bowditch \cite{bowditch-relhyp})-boundary
of $G$, and $\widehat{G} = \GG(G,\HH ) \cup \partial \GG(G,\HH )$ is the compactification of $\partial \GG(G,\HH )$.
The set of  distinct pairs of points on $\partial G$ will be denoted as $\partial^2 G$.

Let $G$ be a (relatively) hyperbolic group. A fixed point on $\partial G$
of a hyperbolic element $g$ of $G$ is called
a {\bf pole} (See \cite{gromov-hypgps} Sec 5.1), and the pair of fixed points $(g_{-\infty}, g_{\infty})$ is called
a {\bf pole-pair}.

\begin{prop} {\bf Pole-pairs Dense} (Gromov \cite{gromov-hypgps} Sec 5.1, p.136)
The collection of pairs $(g_{-\infty}, g_{\infty})$ as $g$ ranges over hyperbolic elements of $G$ is dense in $\partial^2 G$.
\label{poles}
\end{prop}

\subsection{Cannon-Thurston Maps} \label{ct}

Let $(X,{d_X})$ be a proper hyperbolic metric space and $G$ be a (Gromov) hyperbolic group acting freely, properly discontinuously
by isometries on $X$. Let $\Gamma$ be a Cayley graph of $G$ with respect to some finite generating set.
By
adjoining the Gromov boundaries $\partial{X}$ and $\partial{G}$
 to $X$ and $\Gamma$, one obtains their compactifications
$\widehat{X}$ and $\widehat{G}$ respectively.
Choose a basepoint $o\in X$.
Let $ i :\Gamma \rightarrow X$ denote an `inclusion map' mapping $g \in \Gamma$ to $g.o$ and an edge $[a,b]$
of $\Gamma$
to a geodesic  in $X$ joining $a.o, b.o$.

 A {\bf Cannon-Thurston map} $\hat{i}$  from $\widehat{G}$ to
 $\widehat{X}$ is a continuous extension of $i$. The restriction of $\hat{i}$ to $\partial G$ will be denoted by
$\partial i$. The map $\partial i$ induces a relation $\RCT$ on $\partial G$ where $x \sim y$ if
$\partial i (x) = \partial i (y)$ for $x, y \in \partial G$.

\begin{definition}
A {\bf CT leaf} $\lambda_{CT}$ is a bi-infinite geodesic in $\Gamma$
whose end-points are identified by $\partial i$.
\end{definition}

Next, let $G$ be a relatively hyperbolic group, hyperbolic relative to a collection $H_1, \cdots , H_k$
of subgroups.
Let, as before,
$(X,{d_X})$ be a hyperbolic metric space and $G$  act freely, properly discontinuously
by isometries on $X$, such that each conjugate of $H_i$ fixes a unique point of $\partial X$. We thus think of
each conjugate of $H_i$ acting on $X$ as a group of parabolic isometries.
Let, as in Section \ref{ct},
$ i :\Gamma \rightarrow X$ denote the `inclusion map' mapping $g \in \Gamma$ to $g.o$ for some basepoint
$o \in X$. Then $i$ induces a map $i^h:  \GG(G,\HH ) \rightarrow X$.

A {\bf Cannon-Thurston map} $\hat{i}$  from $\widehat{G}$ to
 $\widehat{X}$ is a continuous extension of $i^h$. As before,
 the restriction of $\hat{i}$ to $\partial G$ will be denoted by
$\partial i$ and the induced relation on $\partial G$ by $\RCT$.

Also, the set of  {\it distinct} pairs of points identified by
$\partial i$ will be denoted as $\RCTT$, which is a subset
of $\partial^2 (G)$.

\begin{lemma} $i: \Gamma \rightarrow X$ (or $i^h:  \GG(G,\HH ) \rightarrow X$)
is a quasi-isometric embedding if and only if $\RCTT = \emptyset$. \label{qc} \end{lemma}

\noindent {\bf Proof:} If $i$ (or $i^h$) is a quasi-isometric embedding then $\partial i
: \partial G \rightarrow \partial X$
is a homeomorphic embedding (Theorem 7.2H p. 189 of \cite{gromov-hypgps}). Hence  $\RCTT = \emptyset$.

Conversely, if $i$ (or $i^h$) is {\it not} a quasi-isometric embedding then there exist
$a_n, b_n \in i(\Gamma )$ (or $i^h(  \GG(G,\HH ))$) such that the geodesic $[a_n, b_n]_\Gamma \subset \Gamma$
passes through $1 \in \Gamma$, but the geodesic $[a_n, b_n]_X \subset X$ joining $i(a_n), i(b_n) $
in $X$ lies outside $B_n(i(1))$, the $n$-ball about $i(1)$ in $X$. Assume (after subsequencing)
that $a_n \rightarrow a_\infty
\in \partial G$ and $b_n \rightarrow b_\infty
\in \partial G$. It follows that $\partial i (a_\infty ) = \partial i (b_\infty )$.
Hence $a_\infty \sim b_\infty $ and $\RCTT \neq \emptyset$. $\Box$

\begin{definition} $\RCT$ will be called {\bf trivial} if $x \sim y$ for all $x, y \in \partial G$, and
{quasiconvex} if $\RCTT = \emptyset$. \end{definition}


\subsection{Properties of the $CT$ relation}

We assume in this subsection that $G$ is a (relatively)
hyperbolic group acting freely, properly discontinuously
by isometries on
$(X,{d_X})$, a hyperbolic metric space. Further we  assume that
a  Cannon-Thurston map $\hat{i}$  from $\widehat{G}$ to
 $\widehat{X}$ exists.
The restriction of $\hat{i}$ to $\partial G$ is denoted by
$\partial i$. Recall that $\partial i$ induces a relation $\RCT$ on $\partial G$ where $x \sim y$ if
$\partial i (x) = \partial i (y)$ for $x, y \in \partial G$. Note that $\RCT$ is $G$-invariant.
We investigate the basic properties of $\RCT$ now. Assume that $\RCT$ is non-trivial.

The following Lemma is a direct consequence of the continuity of $\hat i$.

\begin{lemma} \label{closed}
$\RCT$  is a closed relation on $\partial G$, i.e. if $x_n \sim y_n$
for sequences $\{ x_n \} , \{ y_n\} \in \partial G$ and $x_n \rightarrow x \in \partial G$,
$y_n \rightarrow y \in \partial G$, then $x \sim y \in \partial G$. Further,
$\RCTT$ is a closed subset
of $\partial^2 (G)$.\end{lemma}

\begin{definition} $G$ is said to act on $X$ without accidental parabolics   if for every
hyperbolic element $g$ of $G$, its fixed points on $\partial X$ are distinct. \end{definition}

The next Lemma says that poles cannot lie in $\RCT$.

\begin{lemma} Suppose $G$  acts on $X$ without accidental parabolics
If $(x,y) \in \RCT$ and $x \neq y$, then $x$ cannot be a pole of $G$. \label{polerct}\end{lemma}

\noindent {\bf Proof:} By $G$-invariance of $\RCT$ we can assume that the $CT$-leaf $(x,y) \subset
\Gamma$ passes through $1 \in \Gamma$. We argue by contradiction.
If $x$ is a pole, there exists a hyperbolic element $g \in G$ such that $g_\infty = x$ and so $g$ acts
on $X$ as a hyperbolic isometry. Then the infinite geodesic ray $[1,x) \subset \Gamma$
is a $C$-quasigeodesic in $X$ for some $C\geq 1$ since $G$ acts
without accidental parabolics. Choosing a sequence of points $x_n \in [1,x) \subset (x,y)$
such that $x_n \rightarrow \infty$, note that $x_n^{-1}.(x,y)$ and $x_n^{-1}.[1,x)$ converge (up to subsequencing)
to the same bi-infinite geodesic $(p,q) \subset \Gamma$. Hence $(p,q) \subset \Gamma$ is a
$C$-quasigeodesic in $X$. In particular $\partial i (p) \neq \partial i (q)$.
But, by Lemma \ref{closed} $(p,q) \in \RCT$, i.e. $\partial i (p) = \partial i (q)$.
This is a contradiction. $\Box$

\medskip

\noindent {\bf Density} \\ The next Lemma proves density of orbits of cosets of $\RCT$ in the
Hausdorff metric.

\begin{lemma} Let $G$ be a (relatively) hyperbolic group acting on a hyperbolic metric space
without accidental parabolics
and admitting a Cannon-Thurston map. Assume that the relation $\RCT$ is non-trivial, i.e. not all points are
in the same equivalence class.
Let $K \subset \RCT$ be a coset (equivalence class) of the relation.
 Let $C_c(\partial G)$ denote the space of closed subsets
of $\partial G$ with the Hausdorff metric. Then for all $x \in \partial G$, the singleton set
$\{ x \}$ is an accumulation point of $\{ g.K: g \in G \}$. \label{dense} \end{lemma}

\noindent {\bf Proof:}  Clearly,
$K \subset \partial G$ is closed. By Proposition \ref{poles} and Lemma \ref{polerct}, $K$ is nowhere dense
in $\partial G$ and, for all $x \in \partial G$, there exist a sequence $\{ h_n \} \subset G$ such that \\
1) the pole pair
$( h_{n, -\infty}, h_{n, \infty} )$ of $h_n$ is a subset of $\partial G \setminus K$\\
2) $h_{n, \infty} \rightarrow x \in \partial G$ as $n \rightarrow \infty$.\\

By taking sufficiently large powers $h_n^{i_n}$ of $h_n$, we can ensure that $h_n^{i_n}(K)$
lies in an $\epsilon_n$ Hausdorff neighborhood of $h_{n, \infty}$, where $\epsilon_n
\rightarrow 0$. Hence $h_n^{i_n}(K) \rightarrow \{ x \} \in C_c( \partial G )$. $\Box$

\begin{rmk} The assumption on absence of accidental parabolics is not too restrictive.
In all cases of interest, subgroups $H$ of $G$ which become parabolic in $X$ will be regarded
as elements of the the collection $H_1, \cdots , H_n$ relative to which  $G$ is hyperbolic.
\label{densermk} \end{rmk}


\subsection{Continuity in Uniform Topology}
Let $\Lambda_G \subset \partial X$ denote the limit set of $G$, i.e. the collection of accumulation
points in $\partial X$ of an(y) orbit of $G$ acting on $X$.
The existence of a Cannon-Thurston map ensures that $\partial i: \partial G \rightarrow \Lambda_G$
is a quotient map, where pre-images of $\partial i$ are given by elements of $\RCT$.
Let $\partial i_c:C_c( \partial G) \rightarrow C_c(\Lambda_G)$ be the map induced by $i_c$
from $\RCT$-saturated compact subsets of $\partial G$ to
compact subsets of $\Lambda_G$ both equipped with the Hausdorff metric.
Suppose that $f$ is a homeomorphism of $\partial G$ that preserves the cosets of $\RCT$, i.e.
for any coset $K \in \RCT$, $f(K)$ is also a coset of $\RCT$. Then the quotient map from $\Lambda_G$ to itself
induced by $f$ will be denoted by $\overline{f}$.

\begin{prop} \label{unif}
Suppose that $\RCT$ is non-trivial.
Let $f_n$ be a sequence of homeomorphisms of $(\partial G, d)$ that preserves the cosets of $\RCT$, where $d$ denotes some visual metric.
Let $\overline{f_n}$ denote the induced homeomorphisms of $\Lambda_G$.
If $f_n \rightarrow id$ in the uniform topology on $Homeo(\partial G)$ then
$\overline{f_n}  \rightarrow id$ in the uniform topology on $Homeo(\Lambda_G)$. Conversely, if $\overline{f_n}  \rightarrow id$ 
in the uniform topology on $Homeo(\Lambda_G)$ then for every pole $ p \in \partial G$, $d(p, f_n(p)) \rightarrow 0$. More generally, if $\{ p \}$
is a coset of $\RCT$, $d(p, f_n(p)) \rightarrow 0$.
\end{prop}

\noindent {\bf Proof:}  The forward direction follows from the definition of quotient maps.
Conversely, suppose $\overline{f_n}  \rightarrow id$ in the uniform topology on $Homeo(\Lambda_G)$
but $f_n$ does not converge to
the identity map in the uniform topology on $Homeo(\partial G)$. Suppose that there exists (up to subsequencing)
an
$\epsilon > 0$
and a pole $x \in \partial G$ (or a point $x$ such that $\{ x \}$
is a coset of $\RCT$) such that for all $n \in \mathbb{N}$ $d(f_n(x) , x) > \epsilon$.

Then there exists $\eta < \frac{\epsilon}{10}$, say, such that $d_H(B_\eta (x), f_n(B_\eta (x)))  > \eta$
where $d_H$ denotes the Hausdorff metric on $C_c(\partial G)$ and $B_\eta (w)$ denotes a closed ball of radius
$\eta$ centered at $w$.

Assume after subsequencing that $f_n(x) \rightarrow y \in \partial G$.
Then $x \sim y$, since $\overline{f_n}  \rightarrow id$. But $x \neq y$.

This contradicts Lemma \ref{polerct}, if $x$ is  a pole (or more generally, if $\{ x \}$
is a coset of $\RCT$). $\Box$


\section{Commensurators of Kleinian Groups}
In this section we prove discreteness of commensurators for finitely generated, discrete, Zariski dense, infinite covolume Kleinian groups.

\subsection{Surface Groups}

It follows from the Scott core Theorem \cite{scott-cc} that any finitely generated
Kleinian group $H$ is the fundamental group of a compact 3 manifold with boundary. Further, a
geometrically finite group $H_{gf}$ can be chosen abstractly isomorphic to $H$ such that the isomorphism
preserves parabolics (see \cite{mms} for relative versions of the Scott core Theorem).
Hence (see \cite{farb-relhyp} for instance) abstractly $G$ is hyperbolic relative to its parabolic
subgroups. The convex core of ${\mathbb{H}}^3/H$ will be denoted by $M$ and the convex core of ${\mathbb{H}}^3/H_{gf}$ will be denoted by $K$, so that the
inclusion $K \subset M$ is a homotopy equivalence.

\begin{theorem} \cite{mahan-split} \cite{mahan-kl} Let $G$ be a finitely generated Kleinian group. Then $G$ admits
a Cannon-Thurston map from the (relative) hyperbolic boundary of $G$ to its limit set. Further, the Cannon-Thurston map from the 
(relative) hyperbolic boundary of $G$ to its limit set identifies precisely the end-points of leaves of the ending laminations.
\label{ct-kl} \end{theorem}

By Corollary \ref{zariskicor}, it suffices to consider the case where all ends of $M$ are degenerate as these are precisely the manifolds
that come from Kleinian groups with empty domain of discontinuity. 
Our study of $Comm(H)$ splits into two cases: \\
\noindent {\bf Case I:} Some component of the boundary $\partial K$  of $K$ is an incompressible geometrically finite surface \\
{\bf Case II:} $K$ is topologically a compression body whose lower boundary is a (possibly empty) collection of tori.\\

This dichotomy stems from the fact that if some boundary component $Q$ of $K$ is incompressible relative to rank one cusps, then
some other  component of $\partial K$   is an incompressible geometrically finite surface unless $K$ is topologically a 
compression body whose lower boundary is a (possibly empty) collection of tori. This follows from standard 3-manifold topology
\cite{hempel-book} by choosing a maximal collection of disjoint homotopically distinct compressing disks with boundary curves on $Q$.
In particular, if some maximal cusp of $K$ is of rank one, we are in Case I.

We deal with Case I in this subsection and indicate the modifications for Case II in the next subsection. 
By passing to a finite index subgroup of the commensurator group
if necessary, we 
can assume without loss of generality that an end corresponding to a strictly type-preserving representation $\rho : \pi_1(S) \rightarrow Isom({\mathbb{H}}^3)$
is "preserved" by $Comm(H)$. $\rho (\pi_1(S))$ is simply the subgroup of $H$ corresponding to an incompressible end of $M$. By Theorem \ref{ct-kl}
we may assume that $Comm(H)$ preserves the relation $\RCT$ corresponding to the Cannon-Thurston map for the end. Let ${\mathcal L}$ be the ending lamination
whose end-points are identified by the Cannon-Thurston map.

\begin{claim}
Let $\overline{f_n} \in Comm(H)$ be a sequence of commensurators converging to the identity in $Isom({\mathbb{H}}^3)$ and let $f_n$ be the induced homeomorphisms
on the (relative) hyperbolic boundary $\partial \pi_1(S) (=S^1)$ of the group $\pi_1(S)$. Then $f_n \rightarrow Id \in Homeo(S^1)$. \label{claim1}
\end{claim}

\noindent {\bf Proof of Claim:} Suppose $f_n$ does not converge
to the Identity map $Id \in Homeo(S^1)$. Then by Proposition \ref{unif} it suffices to show that there exists (up to subsequencing)
a pole $p \in S^1$ and an $\epsilon > 0$ such that $d(f_n(p),p) \geq \epsilon$ for all $n$, where $d$ is the usual visual metric on $S^1$.

Since $f_n$ does not converge
to the identity, we may assume (after passing to a further subsequence if necessary) that there exists an ideal polygon $\Delta$ such that \\
a) The sides of $\Delta$ are leaves of the ending lamination ${\mathcal L}$\\
b) $f_n|_\Delta \neq Id$ for all $n$.

Then there exist ideal adjacent vertices $A, B$ of $\Delta$ such that $f_n(A) \neq A$ and $f_n(B) \neq B$ for all $n$. Let $\Delta_1$ be another polygon
whose sides lie in ${\mathcal L}$ and whose end-points lie in the arc $(AB)$ of the circle $S^1$ that contains no other ideal points of $\Delta$.
Also let $A_1, B_1$ be adjacent vertices of $\Delta_1$ such that $(A_1B_1) \subset (AB)$. See figure below.

\begin{center}

\includegraphics[height=8cm]{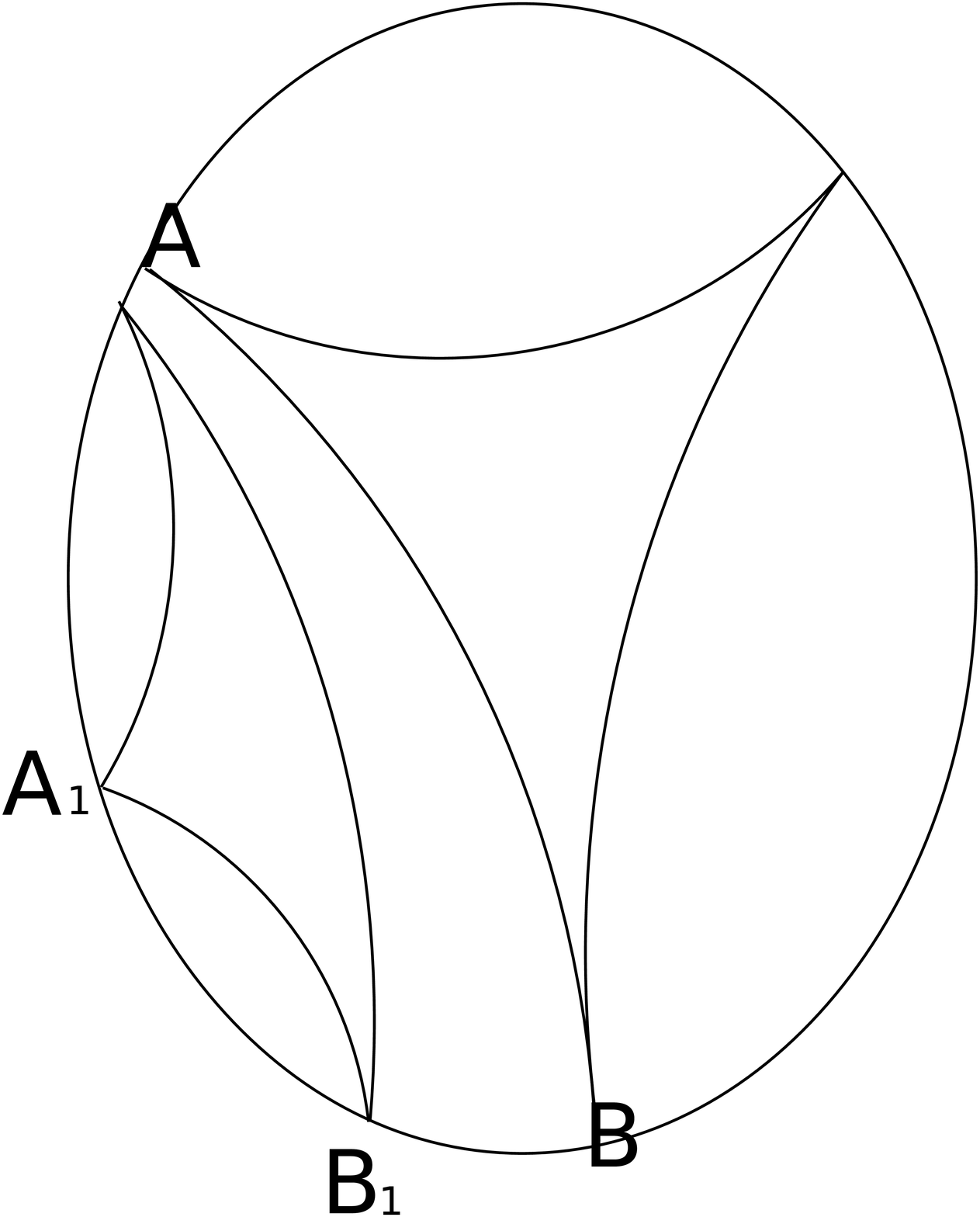}

\end{center}

Since $f_n(A) \neq A$ and $f_n(B) \neq B$ for all $n$
it follows that any pole $p$ in the arc $(A_1B_1)$ is moved at least $\epsilon = min (d(A,A_1), d(A,B_1), d(B,A_1), d(B,B_1))$ by $f_n$ for all $n$.
This proves the Claim. $\Box$

\medskip

We are finally in a position to prove the following. 

\begin{theorem} \label{comm2} Let $G$ be a finitely generated Kleinian group such that
 $G$ is Zariski dense in $PSL_2({\mathbb{C}})$ and ${\mathbb{H}}^3/G$ has infinite volume. Further, suppose that some end of $M = {\mathbb{H}}^3/G$
is incompressible.
Then the commensurator $Comm(G)$ of $G$ is discrete in $PSL_2({\mathbb{C}})$. \end{theorem}

\noindent {\bf Proof:} As above, suppose 
$\overline{f_n} \in Comm(H)$ be a sequence of commensurators converging to the identity in $Isom({\mathbb{H}}^3)$ and let $f_n$ be the induced homeomorphisms
on the (relative) hyperbolic boundary $\partial \pi_1(S) (=S^1)$ of the group $\pi_1(S)$. Then by Claim \ref{claim1} it follows that
for any ideal polygon $\Delta$ with boundary in the ending lamination there exists $N=N(\Delta )$ such that
$f_n$ fixes all the vertices of $\Delta$ for all $n \geq N$. Let $z_\Delta \in S^2_\infty$ be the common image of the end-points of $\Delta$
under the Cannon-Thurston map of Theorem  \ref{ct-kl}. Since the collection of translates of $z_\Delta$ is dense in $S^2$, we can choose
ideal polygons $\Delta_1, \cdots , \Delta_k$ such that the collection of common images $\{ z_1, \cdots , z_k \}$ is Zariski dense in 
$Isom({\mathbb{H}}^3)$. Hence for all $n \geq max_{i=1 \cdots k} \{ N(\Delta_i ) \}$, $\overline{f_n} = Id$. The Theorem follows. $\Box$

\medskip

Theorem \ref{comm2} was proven by Leininger, Long and Reid \cite{llr} under the additional assumption that $G$ is
 non-free and without parabolics.

\subsection{Compressible Core}

It remains to prove the analogue of Theorem \ref{comm2} in Case II, i.e. when the
(relative) core $K$ is topologically a compression body whose lower boundary is a (possibly empty) collection of tori. We describe the
modifications necessary to the proof of Claim \ref{claim1}. Let $S$ be the
upper boundary of $K$ which is a closed surface.
Let $\til K$ be a lift of $K$ to the universal cover $\til M$ of $M$. Then $\til K$ is quasi-isometric to the Cayley graph
of $\pi_1(K)$, which being the fundamental group of a geometrically finite 3-manifold with cusps, is strongly hyperbolic
relative to the maximal parabolic subgroups (see \cite{brahma-pared} or \cite{mahan-split} for instance). Also $\pi_1(K)$ is a free product
of copies of $\mathbb{Z}$ and $\mathbb{Z} \oplus \mathbb{Z}$. In particular, $\pi_1(K)$ has infinitely many ends. 

\medskip

\noindent {\bf Circular order on leaves}\\
The ending lamination 
$\mathcal L$ in this situation
is an element of the Masur domain for $S$ and is well-defined up to Dehn twists about simple closed curves
on $S$ that bound disks in $K$. We first indicate the (relatively simple) argument that works when $K$ is a handlebody, i.e. the collection
of tori on its lower boundary is empty. Here $\pi_1(K)$ is free and is carried by a graph ${\mathcal{G}} \subset S = \partial K$.
Since any commensurator of $G$ induces a homeomorphism of one finite sheeted cover of $K$ to another, it moves $\til K$ to within a bounded Hausdorff distance
of itself in $\til M$. Hence we can isotope the image to $\til K$ by a bounded isotopy. Further, after Dehn twisting about curves in $\partial \til{K}$
that bound disks in $\til K$ we may assume that a lift of ${\mathcal{G}}$ in $\til K$ is isotopic {\bf in $\partial \til{K}$} to its image under
a commensurator followed by an isotopy of $\til K$. Since Dehn twists about curves on $S$ induce homeomorphisms, they induce homeomorphisms of the boundary circle
$S^1$ of $\pi_1(S)$, and hence preserve the circular order on $S^1$. In particular, the nesting of ideal polygons in terms of their ideal vertices
is well-defined. This was {\it the crucial} part of the feature of laminations that we used in the proof of Claim \ref{claim1}.

 We proceed now with the general case, i.e. when the lower boundary of the compression body $K$
might consist of a non-empty collection of tori. We first observe that a leaf of $\mathcal L$ cannot track a horosphere for too long. Suppose not. Then
by $\pi_1(K)$ invariance we can find arbitrarily large segments of leaves in $\til K$ that lie in a bounded neighborhood of a fixed horospherical lift of a
rank 2 cusp in $K$. By taking limits,
we can find a leaf that lies in a bounded neighborhood of a horosphere and hence a leaf of $\mathcal L$ that lies eventually in a
proper subsurface of $S$. This is not possible as $S$ is closed and $\mathcal L$ has no closed leaves. We state this formally below.

\begin{lemma} There exists $D>0$ such that if $l \subset \til{K} $ is a leaf of $\mathcal L$, then the projection $\Pi (l)$ onto any rank 2 cusp
of $\til K$ has diameter bounded by $D$. Further, for every $C>0$, there exists $D_0>0$ such that the intersection $l \cap N_C({\mathcal{H}})$ has diameter
less than $D_0$, where ${\mathcal{H}}$ is any horosphere lift of a rank 2 cusp in $K$. \label{bddcusp} \end{lemma}

Also, by the above discussion, we may assume that every element of the commensurator preserves a natural circular order on the leaves of 
$\mathcal L$ as Dehn twists (the ambiguity up to which $\mathcal L$ is defined) does not change this order. 

Each leaf $l$ of $\mathcal L$ coarsely disconnects 
$\til K$, i.e. for some uniform $\delta > 0$, $\til K \setminus N_\delta(l)$ has infinitely many components. 
Further, by Lemma \ref{bddcusp} we may choose $\delta$ such that each horosphere lift of a rank 2 cusp
lies (coarsely) in exactly one of these components.

Let $\Delta$, $\Delta_1$,
$(AB)$ and $(A_1B_1)$ be as in the proof of Claim \ref{claim1}, where choices are possible due to a well-defined circular order. 
Let  $l_{AB}$ and $l_{A_1B_1}$ be the corresponding leaves of
$\mathcal L$. Then there exists $M_0 \geq 0$ such that for \\
1) a suitable choice of a basepoint $o$ in $\til K$, \\
2)  any component $Q_1$
of $\til K \setminus N_\delta(l_{A_1B_1})$ not containing the  basepoint $o$, and \\
3) any component $Q_2$
of $\til K \setminus N_\delta(l_{CD})$ not containing the  basepoint $o$,
for a side  $l_{CD}$ of $\Delta$ {\it other than} $l_{AB}$, \\
the Gromov inner product $\{u,v\}_o \leq M_0$ for any $u \in Q_1$ and $v \in Q_2$. 

The rest of the argument in Claim \ref{claim1} and Theorem \ref{comm2} go through as before. Hence we have finally,

\begin{theorem} \label{comm2prime} Let $G$ be a finitely generated Kleinian group such that
 $G$ is Zariski dense in $PSL_2({\mathbb{C}})$ and ${\mathbb{H}}^3/G$ has infinite volume. 
Then the commensurator $Comm(G)$ of $G$ is discrete in $PSL_2({\mathbb{C}})$. \end{theorem}

\noindent {\bf Remark:} When $K$ is a handlebody, the circular order is induced on the ends of the free group $\pi_1(K)$ in the lift of the graph
$\mathcal{G} \subset \til{K}$ and we may assume that the leaves of $\mathcal L$ are bi-infinite geodesics on $\mathcal{G} \subset \til{K}$. However, when the 
lower boundary of $K$ consists of a non-empty collection of tori, an alternate argument can be given using the base-points of rank 2 cusps
in the relative hyperbolic boundary of $\pi_1(K)$. The collection of these points on $\partial \pi_1(K)$ are preserved by $Comm(G)$.
Further, as in Claim \ref{claim1} we may assume that there is a fixed ideal polygon $\Delta$ such that for a sequence of commensurators $\overline{f_n}$,
the induced homeomorphisms $f_n$ of $\partial \pi_1(K)$ move $\Delta$. By the circular order on leaves and Lemma \ref{bddcusp} we may assume that
some  rank 2 cusp $\mathcal{H}$ lying in one complementary component of $\til K \setminus N_\delta( \Delta )$ is moved to a different complementary
component. The base-point $p$ of $\mathcal{H}$ in $\partial \pi_1(K)$ is then moved a definite distance under all $f_n$. Then as in the proof of Proposition
\ref{unif} we find that $\overline{f_n} \rightarrow Id$ implies that if $f_n(p) \rightarrow q$ (up to subsequencing)
then the Cannon-Thurston map identifies $p, q$ and hence by Theorem \ref{ct-kl} some leaf of the ending lamination has an ideal end-point
at $p$. This contradicts Lemma \ref{bddcusp}.

\medskip

Using the Thurston-Canary covering theorem \cite{canary-cover}, Theorem \ref{comm2prime} can be strengthened
a la Leininger, Long and Reid \cite{llr} to show that unless $G$ corresponds to a fiber subgroup
of a 3-manifold fibering over the circle, then, in fact, $[Comm(G):G] < \infty$. We discuss this below.

Let $g \in Comm(G)$ be a commensurator of $G$. Then the group $G_1 = <G,g>$ is a finitely generated
discrete Kleinian group. Let $M = {\mathbb{H}}^3/G$ and $M_1 = {\mathbb{H}}^3/G_1$. If $M$ is an
infinite cover of $M$, i.e. if $[G_1:G] = \infty$, then
the Thurston-Canary covering theorem \cite{canary-cover} implies that $G$ is virtually a fiber subgroup
of the finite-volume 3-manifold group $G_1$, i.e. $M_1$ has a finite sheeted cover that fibers
over the circle with fiber subgroup of finite index in $G$.
Else, inductively, for all $g_1, \cdots , g_k \in Comm(G)$, if we denote $G_k = <G,g_1, \cdots , g_k>$, then
$[G_k:G] < \infty$. Let $M_k = {\mathbb{H}}^3/G_1$. Again, by the Thurston-Canary covering theorem \cite{canary-cover},
every end of $M_k$ is covered by an end of $M$. Since $M$ has a finite number of ends, we may assume after
passing to a subsequence, that all the $M_k$ have the same number of ends. Since each end of $M_k$
is homeomorphic to $S \times [0, \infty)$ for some finite volume hyperbolic surface $S$
and since each $S$ can cover only finitely many hyperbolic surfaces, it follows that the
sequence $M_k$ must stabilize. We have thus shown the following.

\begin{theorem} \label{comm3} Let $G$ be a finitely generated Kleinian group such that
 $G$ is Zariski dense in $PSL_2({\mathbb{C}})$ and ${\mathbb{H}}^3/G$ has infinite volume.
Then  $[Comm(G):G] < \infty$ unless $G$ is virtually a fiber subgroup, in which case
$Comm(G)$ is the fundamental group of a virtually fibered finite volume hyperbolic 3-manifold. \end{theorem}

Finally, a note of caution: Just the existence of Cannon-Thurston maps is not sufficient to ensure discreteness of commensurators. A non-uniform arithmetic
lattice in ${\mathbb{H}}^2$ has dense commensurator. However, work of Floyd \cite{Floyd} ensures the existence of Cannon-Thurston maps
in this situation. Thus the proof of Theorem \ref{comm2prime} implicitly  uses the following fact in an essential way:

{\it Any polygon whose sides are leaves of an ending lamination has more than two sides.}

This is exactly what goes wrong for non-uniform lattices in ${\mathbb{H}}^2$ and prevents us from using the separation arguments that go into 
the proof of Theorem \ref{comm2prime}.

\subsection{Examples and Questions}

We  give an example to show that certain naive analogues of Theorem \ref{comm3} cannot be expected in higher dimensions.

\noindent {\bf Example of a group $\Gamma$ whose limit set is a proper subset of $\partial X$, but $[Comm(\Gamma ): \Gamma] = \infty$}\\
Let $M$ be a finite volume hyperbolic 3-manifold fibering over the circle and admitting a totally geodesic embedded
incompressible surface $\Sigma$. Let $\rho : \pi_1(M) \rightarrow Isom({\mathbb{H}}^3)$ be the associated representation.
Let $\rho^{\prime}: \pi_1(M) \rightarrow  Isom({\mathbb{H}}^4) $ be $\rho$ composed with the embedding of $Isom({\mathbb{H}}^3) $
in $Isom({\mathbb{H}}^4) $.
Then a bending deformation of $\rho^{\prime}$, keeping $\rho^{\prime} (\pi_1(\Sigma )$ fixed gives a new discrete faithful representation
$\rho^{\prime \prime}: \pi_1(M) \rightarrow  Isom({\mathbb{H}}^4)  $ with  Zariski dense image
and such that $\rho^{\prime \prime}( \pi_1(M)) = G$ is quasiconvex. Hence $\Lambda_G \neq \partial {\mathbb{H}}^4$. 
Now, let $F$ be the fiber of $M$. Then $\rho^{\prime \prime}( \pi_1(F)) = H$ is normal in $G$ and hence $\Lambda_H = \Lambda_G \neq \partial {\mathbb{H}}^4$.
By Corollary \ref{zariskicor}, $Comm (H)$ is discrete. But $G$ normalizes $H$. Therefore $G \subset Comm (H)$.
Hence $[Comm(H ): H] = \infty$.

\smallskip

\noindent {\bf Absence of examples in dimension 4} \\
We indicate that a naive generalization of Theorem \ref{comm3} cannot exist in dimension 4.
In fact, in \cite{gromov-sq} Gromov comments: \\ "Non-arithmetic $\Gamma$ are especially plentiful for $n = 3$ by Thurston's theory and often
have unexpected features, e.g. some $V = {\mathbb{H}}^3 /\Gamma$ fiber over $S^1$ (which is hard
to imagine ever happening for large $n$). " \\
The next Proposition provides some justification.

\begin{prop}
Let $M$ be a closed aspherical four manifold fibering over the circle with fiber $F$. Then $\pi_1(M)$ cannot be Gromov-hyperbolic. 
In particular $M$ cannot be a rank one locally symmetric space. \label{gromov} \end{prop}

\noindent {\bf Proof:} Observe that $ F \times \mathbb{R}$ covers $M$ and  $\widetilde{ F} \times \mathbb{R}
= \widetilde{M}$ is aspherical. Hence $\pi_i(F) = 0$ for all $i > 1$.  Hence  $F$ is a closed aspherical 3-manifold.
Suppose $\pi_1(M)$ is Gromov-hyperbolic. Then $\pi_1(M)$ cannot contain a copy of $\mathbb{Z} \oplus \mathbb{Z}$. Hence
$F$ is a closed aspherical atoroidal 3-manifold with infinite fundamental group. By Perelman's Geometrization Theorem,
it follows that $F$ admits a hyperbolic structure. Hence $Out (\pi_1 (F))$ is finite. 

But $M = (F \times I) /\phi$ for a diffeomorphism $\phi$ inducing an infinite order element of $Out (\pi_1 (F))$ as otherwise 
$\pi_1(M)$ would contain a copy of $\mathbb{Z} \oplus \mathbb{Z}$.
This is a contradiction. $\Box$

\medskip

Together, Proposition \ref{zariskicor2} and  Theorem \ref{comm2prime} prove discreteness of commensurators for all known examples of finitely
generated, Zariski dense, infinite covolume discrete subgroups of $Isom(X)$ for $X$ a symmetric space of non-compact type. We therefore venture the following

\begin{qn} If $\Gamma$ is a Zariski dense infinite covolume subgroup of $Isom(X)$ for a symmetric space $X$ of non-compact type, is $Comm(\Gamma )$
discrete? \label{qn1} \end{qn}

For finitely generated groups, Question \ref{qn1} subsumes the following (by Proposition \ref{zariskicor2} and  Theorem \ref{comm2prime}).

\begin{qn} If $\Gamma$ is a finitely generated, infinite covolume, 
Zariski dense infinite covolume subgroup of $Isom(X)$ for an irreducible symmetric space $X \neq {\mathbb{H}}^3$ of non-compact type, can the limit
set be all of the Furstenberg boundary? \end{qn}

If we leave the domain of finitely generated groups, then there do exist examples to which  Question \ref{qn1} applies.

\begin{qn} Let $\Gamma_0$ be an infinite index normal subgroup of a rank one (possibly arithmetic) lattice $\Gamma$. 
Is $Comm(\Gamma_0)$ discrete? For instance
$\Gamma_0$ could be the commutator subgroup $[\Gamma , \Gamma ]$ when $H_1(\Gamma , {\mathbb{Q}}) \neq 0$.  \label{qn3} \end{qn}

A particularly simple instance when even the special case of
Question \ref{qn3} seems unknown is when $\Gamma \subset SL_2({\mathbb{R}})$ is a congruence subgroup
of $SL_2({\mathbb{Z}})$.

\medskip

{\bf Acknowledgements:} I would like to thank T. N. Venkataramana, Kingshook Biswas and Pranab Sardar for  helpful email correspondence and discussions.
This work was done in part when the author was visiting Universite Paris-Sud (Orsay) under the ARCUS Indo-French programme. Research partly supported 
by a CEFIPRA project grant.

\bibliography{commens}
\bibliographystyle{alpha}

\end{document}